\documentclass[12pt]{amsart}
\usepackage[T1]{fontenc}
\usepackage{amsmath}
\usepackage{amsthm}
\usepackage{amssymb}
\usepackage{amsfonts}
\usepackage{outlines}
\usepackage{dsfont}
\usepackage{stackengine}
\usepackage{comment}
\usepackage{cancel}
\usepackage{graphicx}
\usepackage{enumerate}
\usepackage{mathrsfs}
\usepackage{tikz-cd}
\usepackage{hyperref}
\usepackage{enumitem}
\usepackage{lipsum}
\usepackage[left = 3.5cm, right = 3.5cm, top = 4 cm, bottom = 4.5 cm]{geometry}

\usepackage[backend=biber, style=alphabetic, sorting=ynt]{biblatex}
\newtheorem{thm}{Theorem} [section]

\newtheorem{lemma}[thm]{Lemma}
\newtheorem{cor}[thm]{Corollary}
\newtheorem{remark}[thm]{Remark}
\newtheorem{prop}[thm]{Proposition}

\title[Decay order bound for mcf near compact singularities]{Decay order bound for mean curvature flow near compact singularities}
\author{Sourav Ghosh}
\date{Nov 19, 2025}

\address{Department of Mathematics, University of Notre Dame, Notre Dame, IN 46556}
\email{sourav.ghosh.1@nd.edu}

\begin{document}

\begin{abstract}
We consider the rescaled flow associated with a mean curvature flow that develops a compact singularity of multiplicity one. We prove that the ``decay order'' of such a rescaled flow is uniformly bounded. As a consequence, we prove a unique continuation result.
\end{abstract}

\maketitle

\section{Introduction}\label{section 1}

A smooth family of hypersurfaces $\{M_t \}_{t \in [0,T)}$ in $\mathbb{R}^{n+1}$ is a mean curvature flow if
$$\frac{\partial}{\partial t} \mathbf x= \mathbf {H}_{M_t} (\mathbf x),$$ 
where $\mathbf {H}_{M_t} (\mathbf x)$ is the mean curvature of $M_t$ at $\mathbf x$. Mean curvature flow is the
gradient flow of area. If the initial hypersurface is smooth and compact, then the flow $M_t$ is known to develop singularities in finite time.

By the works of Huisken \cite{huisken1990asymptotic}, White \cite{white1994partial}, and Ilmanen \cite{ilmanen1995singularities}, it is known that the singularities arising in mean curvature flow can be modeled by self-similar shrinking solutions,
$$\Sigma_t = \sqrt{-t} \Sigma_{-1}, \hspace{0.5 cm} t\in (-\infty,0).$$
To analyze such a singularity, one typically performs a blow-up centered at the point where the singular behavior occurs. This leads to the notion of a tangent flow, introduced in \cite{ilmanen1995singularities}, which generalizes the concept of a tangent cone from the theory of minimal surfaces.

By definition, a tangent flow is obtained as the limit of a sequence of rescalings near the singularity. For example, a tangent flow to the evolving hypersurface $M_t$ at the space-time origin is given by the limit of a sequence of rescaled flows $\frac{1}{\epsilon_i} M_{\epsilon_i^2 t},$ where $\epsilon_i \to 0$. A priori, diﬀerent choices of sequences $\epsilon_i$ may lead to different tangent flows. However, when a tangent flow is modeled on a compact self-shrinker with multiplicity one, Schulze \cite{schulze2014uniqueness} proved that the tangent flow is unique.

For a surface $M$, one can define the Gaussian $L^2$ distance function from $\Sigma$ to $M$ as follows,
$$D_{\Sigma} (M) = \bigg(\int_M d_{\Sigma}^2 \;e^{-\frac{|x|^2}{4}} d \mathcal{H}^n\bigg)^{\frac{1}{2}},$$
where $d_\Sigma$ is the distance function to $\Sigma.$ In practice, a regularized version of this distance function is used, see Section \ref{section 2} for further details. By \cite{schulze2014uniqueness} if the rescaled mean curvature flow \( M_{\tau} \) converges to \(\Sigma\) as \(\tau \to \infty\), then for sufficiently large \(\tau\), \( M_{\tau} \) can be written as the graph of a function \( v(\tau) \) over \(\Sigma\) with 
$$\|v(\tau)\|_{L^2(\Sigma)} \leq C\tau^{\frac{1-2\theta}{\theta}}$$
for some $C>0$ and $\theta \in (0, \frac{1}{2}).$ 

If $M_\tau$ is very close to the shrinker $\Sigma$, one can show that there exists a constant $C_0 > 0$ such that
\[
D_\Sigma(M_{\tau+1}) \leq C_0 \, D_\Sigma(M_\tau),
\]
(see, for instance, Lemma \ref{lemma 2.5} in Section \ref{section 2} and also \cite{ghosh2025cylindrical} for a more general setting). However, this estimate does not rule out the possibility that the decay of the rescaled flow becomes worse as $\tau \to \infty.$ To address this, we aim to establish a uniform bound on the rate at which the rescaled flow approaches the singularity model. In particular, we show that the flow cannot converge to the shrinker too quickly and, in fact, satisfies a doubling-type inequality. This leads to the following result:

\begin{thm}  {\label{thm 1.1}}
    Let $(M_t)_{t \in [-1,0)}$ be a mean curvature flow with a multiplicity one compact singularity modeled by $\Sigma$ at $(0,0)$ and $\{M_\tau \}_{\tau \in [0,\infty)}$ be its associated rescaled mean curvature flow. Then there exists a constant $C > 0$ such that
    $$D_\Sigma(M_{\tau}) \leq CD_\Sigma(M_{\tau+1})$$
    for all $\tau >0.$
\end{thm}

For a rescaled mean curvature flow, one can define the decay order following Sun-Wang-Xue \cite{sun2025passing} of $M_\tau$ by
$$N_{\Sigma}(M_\tau) : = \log \bigg( \frac{D_\Sigma(M_{\tau})}{D_\Sigma(M_{\tau+1})} \bigg).$$
\\
It can be viewed as a discrete parabolic analog of Almgren’s frequency function and doubling constant in the elliptic problems. In essence, this quantity characterizes the asymptotic rate at which the rescaled mean curvature flow converges to the shrinker $\Sigma$. Theorem \ref{thm 1.1} proves that if the singularity is modeled on a compact shrinker, then the decay order $N_{\Sigma}(M_\tau)$ cannot blow up as $\tau \to \infty$.

As a consequence, we obtain the following unique continuation result, which asserts that if the rescaled flow converges to the shrinker at infinite order, then it must coincide with the static solution. This unique continuation property was previously established by Martin-Hagemayer–Sesum \cite{martin2023rate} and Daniels-Holgate-Hershkovits \cite{daniels2025close}.

\begin{cor} {\label{cor 1.2}}
    Let $(M_t)_{t \in (-1,0)}$ be a mean curvature flow with a multiplicity one compact singularity modeled by $\Sigma$ at $(0,0)$ and $\{M_\tau \}_{\tau \in [0,\infty)}$ be its associated rescaled mean curvature flow. Suppose that for all $k > 0$ there exists a constant $C_k$ such that for all $\tau$ we have
    $$D_\Sigma(M_{\tau}) \leq C_k e^{-k\tau},$$
    i.e. the $L^2$ distance from $M_\tau$ to $\Sigma$ vanishes to infinite order as $\tau \to \infty$. Then  $$D_\Sigma(M_{\tau}) =0.$$
\end{cor}

It is an interesting and difficult problem to extend our result to the cylindrical case. In recent work, Sun-Wang-Xue \cite{sun2025regularity} showed that if \( M_\tau \) develops a singularity modeled on a cylinder \( \mathscr{C} \), then as \( \tau \to \infty \), the limit \( \lim_{\tau \to \infty} N_{\mathscr{C}}(M_\tau) \) is either infinite or equal to an eigenvalue of the linearized operator \( L \) (for the definition of \( L \), see Section~\ref{section 2}). Establishing an analogous result in the cylindrical setting would allow one to rule out the possibility that the decay order tends to infinity. We also mention the recent result of Huang-Zhao \cite{huang2025rate}, who showed that if a rescaled mean curvature flow is a global graph over the cylinder, then the distance to the cylinder cannot vanish at infinite order. They also exhibited incomplete graphical solutions that converge to zero at infinite order, with their domains expanding arbitrarily fast.

\subsection{Acknowledgments}
I am grateful to Professor G\'abor Sz\'ekelyhidi for his guidance and support throughout this project. I also thank Professor Nicholas Edelen for many helpful discussions and insights along the way.

\section{Preliminaries}\label{section 2} 

In this section, we will briefly discuss the basic definitions in mean curvature flow that will be useful in the article.

Let $\{M_t\} \subset \mathbb{R}^{n+1}$ be a mean curvature flow and fix a space-time point $(x_0, t_0)$. For any scaling factor $\lambda > 0$, the parabolic rescaling of the flow about $(x_0, t_0)$ with scale $\lambda$ is defined by
\[
M_s^{\lambda} = \lambda \big( M_{t_0 + \lambda^{-2}s} - x_0 \big).
\]
Each rescaled family $\{M_s^{\lambda}\}$ again satisfies the mean curvature flow equation. It allows one to analyze the local behavior of the flow near singularities. For any sequence of such rescaling with $\lambda_i \to \infty,$ there exists a subsequence that converges weakly i.e. as a Brakke flow \cite[Lemma $7.1$]{ilmanen1994elliptic}. Such a flow is called a tangent flow at $(x_0, t_0)$.

We have Huisken’s monotonicity formula \cite{huisken1990asymptotic}:
$$\frac{d}{dt} \int_{M_t} f \rho_{x_0,t_0} d\mathcal{H}^n= \int_{M_t} (\partial_t f - \Delta f)\rho_{x_0,t_0} d\mathcal{H}^n - \int_{M_t} f \bigg| H + \frac{(x -x_0)^{\perp}}{2(t-t_0)}  \bigg|^2 \rho_{x_0,t_0} d\mathcal{H}^n,$$
for $t < t_0$, where $f$ is a function on $M_t$ with polynomial growth (locally uniform in $t$), and
$$\rho_{x_0,t_0} = \frac{1}{\big(4\pi(t_0 - t)\big)^{n/2}} \;\text{exp}\bigg(-\frac{|x-x_0|^2}{4(t_0-t)}  \bigg)$$
is the backward heat kernel. 

The Gaussian weight is 
$$\rho = \frac{1}{(4\pi)^{n/2}} e^{-\frac{|x|^2}{4}}.$$ By $L^2, W^{k,2}$ we denote the weighted Sobolev spaces with respect to $\rho$. 
The Gaussian area functional is 
$$F(M) = \int_{M} \rho \;d\mathcal{H}^n.$$
We define the excess of $M$ to be
$$\mathcal{A}(M) = F(M) - F(\Sigma).$$
Let $(M_t)_{0\leq t < T}$ be a solution to mean curvature flow and let $x_0 \in \mathbb{R}^{n+1}.$  We consider the rescaled flow $M_\tau$ at the space-time point $(x_0,T)$, defined by
$$M_\tau(x) = \frac{1}{\sqrt{T-t}} \big(M_t(x) - x_0 \big), \hspace{0.5cm} \tau(t) = -\log(T-t).$$
Under this rescaling, the flow evolves with normal velocity
$$H + \frac{1}{2}x^\perp.$$ 
The tangent flows can then be analyzed as limits of sequences $M_{\tau_i}$ as $\tau_i \to \infty.$ 
Note that by monotonocity formula for rescaled flow we have,
$$\frac{d}{d\tau} \int_{M_\tau} \rho\;d\mathcal{H}^n = -\int_{M_\tau}  \Big| H + \frac{1}{2}x^\perp \Big|^2 \rho \;d\mathcal{H}^n.$$
Therefore $\mathcal{A}$ is monotonically decreasing along the rescaled flow.

We recall the following lemma from Colding-Minicozzi \cite[Lemma $A.44$] {colding2015uniqueness},

\begin{lemma} {\label{lemma 2.1}}
    For sufficiently small $u$, the graphs $\Sigma_u$ over the surface $\Sigma$ flow by rescaled MCF if and only if $u$ satisfies 
    $$\partial_t u(p,t) = w _u(p) \big(-H_u + \frac{1}{2} \eta_u(p) \big) =: Mu,$$
\end{lemma}

\begin{remark} {\label{remark 2.2}}
For small $u$, the operator $M$ can be decompsoed as  
$$Mu = Lu + Qu,$$ 
where $Lu$ is the linearization of $Mu$ at $u=0$, given by 
$$Lu = \frac{d}{dr} \bigg|_{r=0} M(ru) = \Delta u- \frac{1}{2} x. \nabla u  + |A|^2u  + \frac{1}{2}u,$$
where $A$ denotes the second fundamental form of the underlying hypersurface. The nonlinear term $Qu$ satisfies
$$|Qu| \leq C(|u| + \nabla u|)^2 + C (|u| + |\nabla u|) |\mathrm{Hess}_u|.$$
This decomposition and estimate follow from the calculations in \cite[Appendix A and Lemma $4.10$]{colding2015uniqueness}.
\end{remark}


We say that $\lambda$ is an eigenvalue of operator $L$ if $Lf + \lambda f = 0$ has a nonzero $L^2$ solution. We can decompose \( L^2(\Sigma) \) using an orthonormal basis \( \{ \varphi_i \}_{i=1}^\infty \) with corresponding eigenvalues \( \mu_1 \leq \mu_2 \leq \cdots \to \infty \), such that
\begin{equation}
    L_\Sigma \varphi_i = -\mu_i \varphi_i, \tag{2.2}
\end{equation}
and \( \langle \varphi_i, \varphi_j \rangle = \delta_{ij} \) where the $L^2$-inner product is defined as 
$$\langle u,v \rangle = \int_{\Sigma} u(x) v(x) e^{-\frac{|x|^2}{4}}.$$ 
Using this orthonormal basis, any function $u \in L^2(\Sigma)$ can be written as
\[
    u(x, t) = \sum_{i=1}^{\infty} a_i(t) \varphi_i(x).
\]


\noindent Let $u$ be a solution to the linearized drift heat equation, 
$$(\partial_s - L)u = 0.$$
We consider the $L^2$ norm of $u$, 
$$I(u,s) = \bigg( \int_\Sigma u^2 e^{-\frac{|x|^2}{4}} \bigg)^{\frac{1}{2}}.$$  
It is well known that $\log I(s)$ is convex, and can only be linear when $u = e^{cs} u_0$ for an eigenfunction $u_0$ of $L$ with eigenvalue $c$ (Colding-Minicozzi \cite{colding2020parabolic}). 
In particular for the operators of the form $L_\phi + \lambda$, where 
$$L_\phi(u) = e^\phi\; \text{div} (e^{-\phi} \nabla u),$$
the log-convexity result holds.

The following $3$-annulus type lemma for solutions of our linearized equation has been proved in Lotay–Schulze–Székelyhidi \cite[Lemma $6.1$] {lotay2022neck} (see also \cite[Lemma $3.3$]{simon1984isolated}, \cite[Lemma $5.7$]{ghosh2025cylindrical}). This result is stronger than the log-convexity argument of \cite{colding2020parabolic}.
\begin{lemma} {\label{lemma 2.3}}
    Let $u$ be a solution to the drift heat equation. Then the following hold: 
    \begin{enumerate} [label=\normalfont(\roman*)]
        \item There exists a $\delta>0$ such that for every $\delta_1 < \delta$, one can find a $\delta_2 > \delta_1$ depending on $\delta_1$, with the following property: If $I(u,t+1) \geq e^{\delta_1}  I(u,t)$, then $I(u,t+2) \geq e^{\delta_2}  I(u,t+1).$ 
        \item Suppose $u$ does not have any component with zero eigenvalue. Then there exists a small $\delta >0$ such that $I(u,t+1) \geq e^{\delta}  I(u,t),$ or $I(u,t-1) \geq e^{\delta}I(u,t).$
    \end{enumerate}
\end{lemma} 

We need a quantitative version of this statement. This is similar to Proposition $2.7$ in Sz\'ekelyhidi \cite{szekelyhidi2021minimal}.

\begin{prop} {\label{prop 2.4}}
Let $C_0$ be the constant from Lemma \ref{lemma 2.5}, and let $u$ be a solution to the drift heat equation. Choose any \( L_0 \in \left(0, \tfrac{1}{2} \right) \). Then there exist a constant \( B > 0 \) and a choice of \( L \in \left[L_0/2, L_0 \right] \), depending only on the shrinker $\Sigma$, such that the following holds. Suppose that 
$$I(u,0) \leq 2C_0,$$
and 
$$I(u,2L) \leq \frac{1}{2C_0}.$$
Then we have
$$I(u,L) \leq 1- L^B. $$
\end{prop}

\begin{proof}
By the $L^2$-spectral decomposition as above, we can write
$$u= \sum_j a_j e^{-\lambda_j t}\phi_j.$$
where \( \lambda_1 < \lambda_2 < \cdots \) denote the eigenvalues of the operator \( -L \), and \( \varphi_k \) the corresponding eigenfunctions. Then, by the Plancherel identity, we have
$$I(u,t)^2 = \sum_j a_j^2 e^{-2\lambda_j t}.$$
Then by the Weyl law there exist \( C_1, m > 0 \) such that for any \( N > 1 \), the number of eigenvalues \( \lambda_j \) with \( \lambda_j \leq N \) is at most \( C_1 N^m \). Therefore, given $L_0 \in (0, \frac{1}{2})$, we can choose $L \in [L_0/2, L_0]$ such that
\[
\min_j \left|\lambda_j - L^{-1} \ln(2C_0) \right| \geq \frac{C_1^{-1} L^{-1}}{(L^{-1})^m} = C_1^{-1} L^{m-1}.
\]
We can write the assumptions as
$$\sum_j a_j^2 \leq 4C_0^2,$$
$$\sum_j a_j^2 e^{-4\lambda_jL} \leq \frac{1}{4C_0^2}.$$
It follows that,
$$\sum_j a_j^2 e^{-\ln 4C_0^2} \leq 1,$$
$$\sum_j a_j^2 e^{-4\lambda_jL + \ln 4C_0^2} \leq 1.$$
So we can add two inequalities to get 
$$\sum_j a_j^2 e^{-2\lambda_jL} \cosh (2\lambda_jL - \ln 4C_0^2) \leq 1.$$
We have
$$|2\lambda_jL- \ln (4C_0^2)| = 2L|\lambda_j - L^{-1} \ln (2C_0)| \geq C_1^{-1} L^m.$$
This implies that
$$\cosh (2\lambda_jL - \ln (4C_0^2)) \geq 1 + C_1^{-1} L^m,$$
increasing $C_1$ if necessary and since $L \leq \frac{1}{2},$ so we can choose a constant $B$ large enough depending on $m$ such that 
$$\frac{1}{1+C_1^{-1} L^m} < (1-L^B)^2,$$
which implies,
$$\sum_j a_j^2 e^{-2\lambda_jL} \leq (1-L^B)^2.$$
Hence, the proposition follows.
\end{proof}
Now consider the parabolic equation 
$$(\partial_t - L) u = f, \hspace{0.4 cm} u(x,0) = 0.$$
We will show the operator $\partial_t - L$  has a right inverse. This is just an adaptation of Galerkin’s method from linear parabolic parabolic differential equation.

If we start with a function $f$ such that $f(.,t) \in L^2$ for all $t.$ Then $f$ can be expanded in terms of an orthonormal basis $\{\phi_i\}_{i \geq 0}$ as
$$f(x,t) = a_0(t)\phi_0(x) + \sum_{i=1}^\infty a_i(t) \phi_i(x).$$ 
Similarly, let 
$$u(x,t) = b_0(t)\phi_0(x) + \sum_{i=1}^\infty b_i(t) \phi_i(x),$$ such that
$$b_0(t) = \int_0^t a_0(\tau) d\tau,$$
and for $i \geq 1,$
$$b_i(t) = e^{-\lambda_i t} \int_0^t e^{\lambda_i \tau} a_i(\tau) d\tau.$$
Then $u$ will satisfy the equation $(\partial_t - L)u = f$ in the weak sense. One begins by constructing weak $L^2$ solutions (Evans \cite{evans2022partial}, $7.1.2$, Theorems $3, 4$) and then improves their regularity to obtain strong solutions \cite[$7.1.3$, Theorem $5$] {evans2022partial}. Now we will assume that $0 < t \leq 1.$
Note that, 
$$|b_0(t)| \leq t \sup_{\tau \in [0,t]} \|a_0(\tau)\|,$$
and for $i \geq 1,$
$$|b_i(t)| \leq \frac{1}{\lambda_i}|1-e^{-\lambda_it}| \sup_{\tau \in [0,t]} \|a_i(\tau)\| $$
As we assumed $t \leq 1$ and $\lambda_i \to \infty,$ therefore 
$$|b_i(t)| \leq Ct \sup_{\tau \in [0, t]} \|a_i(\tau)\| $$
Therefore, we can find a right inverse of the operator for $t \leq 1$ such that 
$$\|u(.,t)\|_{L^2} \leq Ct \sup_{\tau \in [0,t]} \|f(.,t)\|_{L^2}.$$
Now we will define our distance function $D_\Sigma(M)$ as follows. Fix an integer $k$ large (which will be chosen in Proposition \ref{prop 3.1}) and a constant $c_0>0$ (which will be chosen in Lemma \ref{lemma 2.5}). Let $M$ be a hypersurface that can be written as a graph over $\Sigma$ of a function $u$ with $\|u\|_{C^{k,\alpha}} < c_0,$ then  we define 
$$D_\Sigma(M) = \bigg(\int_{M} |u|^2 e^{-\frac{|x|^2}{4}} \bigg)^{\frac{1}{2}}. $$
If $M$ is not such a small graph over $\Sigma$, then we let $D_\Sigma(M) = \infty$. 

Note that $|u|(x) = d_{\Sigma}(x+ u(x))$, and under the parametrization of $M$ using $u$, the volume forms of $M$ and $\Sigma$ are uniformly equivalent. Therefore, the integrals of $d^2_{\Sigma} e^{-|x|^2/4}$ over $M$ and of $|u|^2 e^{-|x|^2 /4}$ over $\Sigma$ are uniformly equivalent.
The $L^2$-distance to $\Sigma$ has the following properties.
  
\begin{lemma} {\label{lemma 2.5}}
    There exist constants $\epsilon, C_0 >0$ such that if $D_{\Sigma}(M_\tau), \mathcal{A}(M_\tau) - \mathcal{A}(M_{\tau+1}) < \epsilon,$ then for any $s \in [0,1]$
    $$D_\Sigma(M_{\tau + s}) \leq C_0\, D_\Sigma(M_\tau).$$
\end{lemma}

\begin{proof}
    First note that if $\epsilon>0$ sufficiently small, then by definition of $D_\Sigma(M_\tau),$ we can write $ M_\tau$ as a graph over $\Sigma$ of a function with $C^{k,\alpha}$ norm bounded by $c_0.$ Using the pseudolocality argument (see, for instance, Ilmanen-Neves-Schulze \cite{ilmanen2019short}) we then we can conclude that if $c_0$ small, $M_{\tau+s}$ can also be written as a graph of a function $v$ satisfying $\|v(.,s)\|_{C^{k,\alpha}} < c_1$ for some constant $c_1$.
    
    Next, we show that if $D_{\Sigma}(M_\tau), \mathcal{A}(M_\tau) - \mathcal{A}(M_{\tau+1}) <\epsilon$ and $\epsilon$ is sufficiently small, then we can write $M_{\tau+s}$ as a graph over $\Sigma$ of a function $v(.,s)$ such that $\|v(.,s)\|_{C^{k,\alpha}} < c_0$ for all $s \in [0,1].$ The proof proceeds by a standard contradiction argument. suppose this is not the case. Then there exists a sequence of flows $M^i$ for which $D_{\Sigma}(M^i_\tau), \mathcal{A}(M^i_\tau) - \mathcal{A}(M^i_{\tau+1}) <\epsilon_i \to 0$ such that $M^i$ is a graph over $\Sigma$ of a function $v^i$ whose $C^{k, \alpha}$ norm is bounded by $c_0$ for some $s \in [0,1].$ 
    
    By the compactness theorem for Brakke flows (Ilmanen \cite{ilmanen1994elliptic}, Lemma $7.1$) there exists a subsequence of $\{M^i_{\tau+s}\}_{s \in [0,1]}$ that converges to a limiting Brakke flow. Since $\mathcal{A}(M^i_\tau) - \mathcal{A}(M^i_{\tau+1}) \to 0$ the limit must be a self-shrinker. Moreover, by our assumptions, $\{M^i_{\tau+s}\}_{s \in [0,1]}$ can be written as a $C^{k,\alpha}$ graph over $\Sigma$ of the function $v^i(.,s)$ with $v^i(.,0) \to 0.$ Therefore, applying Allard's regularity theorem \cite{allard1972first}, we obtain $\|v^i(.,s)\|_{C^{k,\alpha}} \to 0$ for all $s$ which contradicts the assumption that for each $i,$ there exists some $s \in [0,1]$ for which $\|v^i(.,s)\|_{C^{k,\alpha}} > c_0$. This completes the argument. 

    Now by Lemma \ref{lemma 2.1} Remark \ref{remark 2.2} , we have
    $$|\partial_s v_s - \Delta v_s + \frac{1}{2} x. \nabla v_s| \leq C \big(|v_s| + |\nabla v_s| \big).$$
    To absorb the gradient term, we will compute the evolution of $|v|^{3/2}$ and obtain a differential inequality of the form
    $$\partial_\tau |v_s|^{3/2} - \Delta |v_s|^{3/2} + \frac{1}{2} x. \nabla |v_s|^{3/2} \leq C |v_s|^{3/2},$$
    for possibly a larger constant $C$. 
    
    Define the space-time function,
    $$f(x,s) = e^{-Cs} |v_s|^{3/2}.$$
    Then $f$ is a subsolution of the drift heat equation on the time interval $[0,1].$ Let 
    $$d := D_\Sigma(M_\tau).$$ 
    By Ecker's log-Sobolev inequality \cite{ecker2000logarithmic} we obtain for any $s \in [0,1],$  
    $$\int_{\Sigma} f^\frac{4q}{3} e^{-\frac{|x|^2}{4}} \leq C d^{2q},$$
    where $q=q(s) \geq 1$ depends on $s$ and satisfies $q(s) \to 1$ as $s \to 0.$ This estimate implies 
    $$\bigg(\int_{M_{\tau+s}} d^{2q(s)} e^{-\frac{|x|^2}{4}} \bigg)^{\frac{1}{q(s)}} \leq D_\Sigma(M_\tau).$$
    Applying Hölder’s inequality, we conclude that 
    $$D_\Sigma(M_{\tau + s}) \leq C_0\, D_\Sigma(M_\tau),$$
    for some constant $C_0 >0$ independent of $s.$
\end{proof}
\section{Proof of the main theorem}\label{section 3}

In this section, we begin by proving a preliminary result that will be essential in the proof of Theorem \ref{thm 1.1}. We will first prove a result similar to the monotonicity of the frequency function in classical proofs of a strong unique continuation. This is analogous to \cite[Proposition $4.4$]{szekelyhidi2021minimal}.

\begin{prop} {\label{prop 3.1}}
    Let $C_0$ and $\epsilon$ be the constant from Lemma \ref{lemma 2.5}. Suppose that $D_{\Sigma}(M_\tau), \mathcal{A}(M_\tau) - \mathcal{A}(M_{\tau+1}) < \epsilon$. Then there exists an $A>0$ such that for any $L_0 \in (0, 1/2)$, one can choose an $L \in [L_0/2, L_0]$ with the following property. Suppose that 
    $$\frac{1}{2C_0} D_\Sigma(M_{\tau}) \leq D_\Sigma(M_{\tau+L}) \leq L_0^A.$$
    Then
    $$\frac{1}{2C_0} D_\Sigma(M_{\tau + L}) \leq D_\Sigma(M_{\tau + 2L}).$$
\end{prop}

\begin{proof}
    We will prove the result by contradiction. Let us denote 
    $$a:= D_\Sigma(M_{\tau+L}).$$
    Suppose that
    $$D_\Sigma(M_{\tau}) \leq 2aC_0,$$
    and 
    $$D_\Sigma(M_{\tau+2L}) \leq \frac{a}{2C_0}.$$
    We will show that there exists a constant $A>0$ such that if $a < L_0^A$, then for a suitable choice of $L \in [L_0/2, L_0],$ we have 
    $$D_\Sigma(M_{\tau+L}) < a,$$ 
    which contradicts the definition of $a.$ By Lemma \ref{lemma 2.5}, we can write $M_{\tau+s}$ as a graph over $\Sigma$ of a function $v(s)$ for $s \in [0,1]$ with $\|v(s)\|_{C^{k,\alpha}} \leq c_0.$ Define, 
    $$Q(v):= (\partial_s - L) v.$$
    By the interpolation inequality (see, for example, \cite{colding2015uniqueness}), we have
    $$|\nabla v| \leq C a^{\frac{k-1}{k+n+1}}, \hspace{0.5 cm} |\nabla^2 v| \leq C a^{\frac{k-2}{k+n+1}}.$$
    Hence, by Remark \ref{remark 2.2}, we can choose $k$ sufficiently large so that $\|Qv\|_{L^2} \leq C a^{1+\gamma}$ for some $\gamma >0$. We can solve the equation
    $$(\partial_s - L) w = Qv, \hspace{0.5 cm} w(0) = 0.$$
    Note that 
    $$I(w,L), I(w,2L) \leq CLa^{1+\gamma}.$$
    Define 
    $$V:= v - w$$
    Then $V$ satisfies
    $$(\partial_s - L)V = 0.$$
    By construction,
    \begin{align*}
        I(V,0) &\leq 2aC_0 + Ca^{1+\gamma}, \\
        I(V,2L) &\leq \frac{a}{2C_0} + CLa^{1+\gamma}.
    \end{align*}
    Set 
    $$C_a= a + 2C_0CLa^{1+\gamma}, \hspace{0.5 cm} \widehat{V} = C_a^{-1}V.$$
    Then
    \begin{align*}
        I(\widehat{V},0) &\leq 2C_0, \\
        I(\widehat{V},2L) &\leq \frac{1}{2C_0}.
    \end{align*}
    This is exactly the setting of Proposition \ref{prop 2.4}. Therefore, for an appropriate choice of $L \in [L_0/2, L_0]$, 
    $$I(\widehat{V},L) \leq 1 - L^B.$$
    for some $B>0.$ Rescaling back gives
    $$I(V,L) \leq (a+2CC_0La^{1+\gamma})(1 - L^B),$$
    which implies
    $$D_\Sigma(M_{\tau+L}) \leq (a+2CC_0La^{1+\gamma})(1 - L^B) + CLa^{1+\gamma}.$$
    Thus, for some $C_1> 0,$ 
    $$D_\Sigma(M_{\tau+L}) \leq a - aL^B + C_1La^{1+\gamma}.$$
    Since 
    $$a = D_\Sigma(M_{\tau+L}) \leq C_0\epsilon,$$
    we can choose $0< \gamma_0 < \gamma$ and $\epsilon$ small enough to absorb the constant $C_1$, ensuring that
    $$a^{\gamma_0} > C_1a^{\gamma}.$$ 
    Next, we choose 
    $$A = \frac{B-1}{\gamma_0}.$$ 
    Then if $a \leq L^A$, we obtain
    $$D_\Sigma(M_{\tau+L}) < a,$$ 
    which contradicts the definition of $a$ at the beginning of the proof.
\end{proof}

Now we can prove our main Theorem $\ref{thm 1.1}$.

\begin{proof}[Proof of Theorem $\ref{thm 1.1}$] 
Since $\Sigma$ is the unique tangent flow of $M_\tau,$ we can first choose $T_0$ large enough such that for $T \geq T_0, D_{\Sigma}(M_\tau), \mathcal{A}(M_\tau) - \mathcal{A}(M_{\tau+1}) < \epsilon,$ where $\epsilon$ is from Lemma $\ref{lemma 2.5}$. Now, let $L_0 < \frac{1}{2}.$ By increasing $T_0$ if necessary, we can ensure $D_{\Sigma}(M_\tau) < L_0^A,$ where $A$ is from Proposition \ref{prop 3.1}. Let $C_0$ be the constant from Lemma \ref{lemma 2.5}.

Suppose there exists a time $T \geq T_0$ such that for all $L \in [L_0/2, L_0]$, the following inequality holds: 
$$D_\Sigma(M_{T}) \leq 2C_0D_\Sigma(M_{T+L}).$$ 
Then, based on our assumptions, Proposition \ref{prop 3.1} implies that for all $m \in \mathbb{N}$, we have
$$D_\Sigma(M_{T + mL}) \leq 2C_0 D_\Sigma(M_{T + (m+1)L}).$$
Now, let $\tau \geq T.$ We can choose a $m \in \mathbb{N}$ depending on $L_0$ and $N \in \mathbb{N}$ depending on $\tau,$ such that 
$$T+NL_0 \leq \tau \leq T + (N+1)L_0\leq T+ (N + m-1)L \leq \tau+1 \leq T+ (N + m)L.$$
Consequently, 
\begin{align*}
    D_\Sigma(M_\tau) &\leq C_0 D_\Sigma(M_{T+NL_0}) \\
    &\leq 2C_0.C_0 D_\Sigma(M_{T+ (N+1)L_0}) \\
    &\leq (2C_0)^m C_0 D_\Sigma(M_{T+ (N + m)L_0}) \\
    &\leq (2C_0)^m C_0^2 D_\Sigma(M_{\tau+1}).
\end{align*}
Thus we obtain
$$D_\Sigma(M_\tau) \leq A_0 D_\Sigma(M_\tau),$$
for some constant $A_0 >0$, for all $\tau \geq T_0.$

Therefore, we can assume that for all $T > T_0$, there exists some $L \in [L_0/2, L_0]$ such that
$$D_\Sigma(M_{T}) \geq 2C_0 D_\Sigma(M_{T+L}).$$
Since, $L < \frac{1}{2}$, by Lemma \ref{lemma 2.5}, we have
$$D_\Sigma(M_{T+L}) \geq C_0^{-1} D_\Sigma(M_{T+L_0}),$$
and combining this with the previous inequality, we obtain
$$D_\Sigma(M_T) \geq 2D_\Sigma(M_{T+L_0}).$$
Thus, we can bound $D_\Sigma(M_{T+L_0})$ as
$$D_\Sigma(M_{T+L_0}) \leq \frac{1}{2}D_\Sigma(M_T) \leq \frac{1}{2} L_0^A.$$
Now, define 
$$L_1 := \bigg(\frac{1}{2} \bigg)^{\frac{1}{A}} L_0.$$
This gives the estimate
$$D_\Sigma(M_{T+L_0}) \leq L_1^A.$$
We can now repeat the construction by replacing $M_T$ with $M_{T+L_0}$ and $L_0$ with $L_1.$ In particular, if for any $L \in [L_1/2, L_1]$ the condition
$$D_\Sigma(M_{T+L_0}) \leq 2C_0 D_\Sigma(M_{T+L_0+L})$$
holds, then, as before, we obtain the bound
$$D_\Sigma(M_\tau) \leq A_1 D_\Sigma(M_{\tau+1}),$$
for some constant $A_1>0$ (possibly $A_1 > A_0,$ with $A_1$ depends on $L_1$) for all $\tau \geq T+L_0$. 

Otherwise, we will have 
$$D_\Sigma(M_{T+L_0+L_1}) \leq \frac{1}{2} L_1^A.$$
We will iterate this process and it will eventually stop as 
$$\sum_{i=0}^\infty L_i =  \frac{L_0}{1- (\frac{1}{2})^{\frac{1}{A}}} < \infty. $$ \\
so it will imply $D_\Sigma(M_\tau) = 0$ for some large $\tau.$ However, this would mean that $M_\tau$ is the static flow.
\end{proof}
Based on this result, we can prove Corollary \ref{cor 1.2}.

\begin{proof}[Proof of Theorem $\ref{cor 1.2}$] 
Suppose that for all $k > 0$ there exists a constant $C_k$ such that for all $\tau$ we have
$$D_\Sigma(M_{\tau}) \leq C_k e^{-k\tau}.$$
From Theorem \ref{thm 1.1}, we have
$$D_\Sigma(M_0) \leq C D_\Sigma(M_1).$$
By induction, for all $L \in \mathbb{N}$,
$$D_\Sigma(M_0) \leq C^L D_\Sigma(M_L).$$
Therefore using our assumption we can now bound $D_\Sigma(M_0)$ as follows,
$$D_\Sigma(M_0) \leq C^L C_ke^{-kL}.$$
The inequality holds for all $k>0$ and $L>0$. First, choose $k$ large enough so that $Ce^{-k}< \frac{1}{2}$. Then for this choice of $k,$ the right-hand side of the inequality becomes arbitrarily small as $L$ increases. Therefore, we conclude that $D_\Sigma(M_0) = 0.$
\end{proof}

\printbibliography

@article{colding2015uniqueness,
  title={Uniqueness of blowups and {\L}ojasiewicz inequalities},
  author={Colding, Tobias Holck and Minicozzi, William P},
  journal={Annals of Mathematics},
  pages={221--285},
  year={2015},
  publisher={JSTOR}
}

@article{colding2020parabolic,
  title={Parabolic frequency on manifolds},
  author={Colding, Tobias Holck and Minicozzi II, William P},
  journal={arXiv preprint arXiv:2002.11015},
  year={2020}
}

@article{lotay2022neck,
  title={Neck pinches along the Lagrangian mean curvature flow of surfaces},
  author={Lotay, Jason D and Schulze, Felix and Sz{\'e}kelyhidi, G{\'a}bor},
  journal={arXiv preprint arXiv:2208.11054},
  year={2022}
}

@book{ilmanen1994elliptic,
  title={Elliptic regularization and partial regularity for motion by mean curvature},
  author={Ilmanen, Tom},
  volume={520},
  year={1994},
  publisher={American Mathematical Soc.}
}

@article{ilmanen2019short,
  title={On short time existence for the planar network flow},
  author={Ilmanen, Tom and Neves, Andr{\'e} and Schulze, Felix},
  journal={Journal of Differential Geometry},
  volume={111},
  number={1},
  pages={39--89},
  year={2019},
  publisher={Lehigh University}
}

@article{ecker2000logarithmic,
  title={Logarithmic Sobolev inequalities on submanifolds of Euclidean space},
  author={Ecker, Klaus},
  year={2000},
  publisher={Walter de Gruyter GmbH \& Co. KG Berlin, Germany}
}

@article{sun2025regularity,
  title={Regularity of cylindrical singular sets of mean curvature flow},
  author={Sun, Ao and Wang, Zhihan and Xue, Jinxin},
  journal={arXiv preprint arXiv:2509.01707},
  year={2025}
}

@article{sun2025passing,
  title={Passing through nondegenerate singularities in mean curvature flows},
  author={Sun, Ao and Wang, Zhihan and Xue, Jinxin},
  journal={arXiv preprint arXiv:2501.16678},
  year={2025}
}

@article{huisken1990asymptotic,
  title={Asymptotic behavior for singularities of the mean curvature flow},
  author={Huisken, Gerhard},
  journal={Journal of Differential Geometry},
  volume={31},
  number={1},
  pages={285--299},
  year={1990},
  publisher={Lehigh University}
}

@article{ilmanen1995singularities,
  title={Singularities of mean curvature flow of surfaces},
  author={Ilmanen, Tom},
  journal={preprint},
  volume={1},
  pages={23--25},
  year={1995}
}

@article{schulze2014uniqueness,
  title={Uniqueness of compact tangent flows in mean curvature flow},
  author={Schulze, Felix},
  journal={Journal f{\"u}r die reine und angewandte Mathematik (Crelles Journal)},
  volume={2014},
  number={690},
  pages={163--172},
  year={2014},
  publisher={De Gruyter}
}

@article{white1994partial,
  title={Partial regularity of mean-convex hypersurfaces flowing by mean curvature},
  author={White, Brian},
  journal={International Mathematics Research Notices},
  volume={1994},
  number={4},
  pages={185--192},
  year={1994},
  publisher={Citeseer}
}

@book{evans2022partial,
  title={Partial differential equations},
  author={Evans, Lawrence C},
  volume={19},
  year={2022},
  publisher={American Mathematical Society}
}

@article{ghosh2025cylindrical,
  title={Cylindrical tangent flows in mean curvature flow},
  author={Ghosh, Sourav},
  journal={arXiv preprint arXiv:2508.05517},
  year={2025}
}

@article{szekelyhidi2021minimal,
  title={Minimal hypersurfaces with cylindrical tangent cones},
  author={Sz{\'e}kelyhidi, G{\'a}bor},
  journal={arXiv preprint arXiv:2107.14786},
  year={2021}
}

@article{daniels2025close,
  title={How close is too close for singular mean curvature flows?},
  author={Daniels-Holgate, Joshua and Hershkovits, Or},
  journal={arXiv preprint arXiv:2503.11522},
  year={2025}
}

@article{martin2023rate,
  title={On the rate of convergence of the rescaled mean curvature flow},
  author={Martin-Hagemayer, Rory and Sesum, Natasa},
  journal={arXiv preprint arXiv:2302.06742},
  year={2023}
}

@article{huang2025rate,
  title={On the rate of convergence of cylindrical singularity in mean curvature flow},
  author={Huang, Yiqi and Zhao, Xinrui},
  journal={arXiv preprint arXiv:2510.23499},
  year={2025}
}

@article{allard1972first,
  title={On the first variation of a varifold},
  author={Allard, William K},
  journal={Annals of mathematics},
  volume={95},
  number={3},
  pages={417--491},
  year={1972},
  publisher={JSTOR}
}

@incollection{simon1984isolated,
  title={Isolated singularities for extrema of geometric variational problems},
  author={Simon, Leon},
  booktitle={Miniconference on nonlinear analysis},
  volume={8},
  pages={46--51},
  year={1984},
  publisher={Australian National University, Mathematical Sciences Institute}
}

\end{document}